\documentclass[11pt]{amsart}
\usepackage{amsfonts}
\usepackage{amssymb}

\def\la{\lambda}	 		
	 			\def\t{\theta}
\def\f{\phi}

	 \def\C{{\mathbb C}}
\def\D{{\mathbb D}}

\def\T{{\mathbb T}}	 

\def\({\left(}		 \def\){\right)}

\newtheorem{thm}{\sc Theorem}
\newtheorem{cor}{\sc Corollary}

\begin{document}
\title[Taylor coefficients of inner functions and Beurling's theorem]
{Taylor coefficients of inner functions and Beurling's theorem for the shift}
\author[D. Vukoti\'c]{Dragan Vukoti\'c}
\address{Departamento de Matem\'aticas, Universidad Aut\'onoma de
Madrid, 28049 Madrid, Spain}
\email{dragan.vukotic@uam.es}
\urladdr{http://www.uam.es/dragan.vukotic}
\date{30 December, 2014}
\begin{abstract}
The purpose of this expository note is to give a proof of a Schur-type theorem that characterizes the inner functions in terms of their Taylor coefficients. In view of Beurling's theorem, this provides a sequential characterization of the shift-invariant subspaces of $\ell^2$.
\end{abstract}
\maketitle
\subsection*{Introduction}
A formal power series $f(z)=\sum_{n=0}^\infty a_n z^n$ represents an  analytic function in the unit disk $\D$ if and only if $\limsup_{n \to\infty} |a_n|^{1/n}\le 1$. Among such functions, those which are bounded by one (\textit{i.e.\/}, map $\D$ into itself) can be characterized in terms of their Taylor coefficients $(a_n)_n$ by a classical theorem due to Schur \cite{S}. This criterion is expressed by infinitely many conditions; see \cite[p.~40, 180]{G}, \cite{K} or \cite[Theorem~IV.25]{T}. 
\par 
An important subclass of functions that map $\D$ into itself is formed by \textit{inner functions\/}. These are functions whose boundary values have modulus one almost everywhere on the unit circle. Basic facts about inner functions can be found in \cite{D} or \cite{G}. Taylor coefficients of inner functions have been studied in \cite{AJ} and \cite{NS}. By using an approach analogous to that of Kortram \cite{K}, one can obtain a complete characterization of inner functions in terms of their Taylor coefficients; like in Schur's criterion, it is also given in terms of infinitely many conditions. The result appear to be essentially known to the experts and somewhat different proofs are also possible.
\par
The question of describing the lattice of all closed non-trivial invariant subspaces of a given operator is fundamental and such descriptions are available for some classical operators. For the \textit{shift operator\/} which transforms $f(z)$ into $z f(z)$, acting on the Hardy space $H^2$, Beurling's celebrated theorem (\cite{B}, \cite[Chapter~7]{D}, \cite[Theorem~II.7.1]{G}) states that all of its invariant subspaces are of the form $\f H^2$ for some inner function $\f$. The importance of Beurling's theorem stems for the fact that it gives a function-theoretic meaning to each invariant subspace. The space $H^2$ is usually identified with the complex space $\ell^2$ of all square-summable sequences and the shift with the operator $S$ on $\ell^2$ which acts as follows: $S(w_0,w_1,w_2,\ldots)= (0,w_1,w_2,\ldots)$. The characterization of inner functions allows for a sequential description of invariant subspaces. Perhaps it is not very practical to check but nonetheless one can write down such a statement.
\par
\subsection*{Preliminaries}
Let $H^\infty$ denote the Banach space of bounded analytic functions in $\D$ equipped with the usual norm $\|f\|_\infty = \sup_{z\in\D} |f(z)|$. The Hardy space $H^2$ of the disk is the Hilbert space of all functions $f$ analytic in $\D$ for which
$$
 \|f\|_2 = \lim_{r\to 1} \( \int_0^{2\pi}|f(r e^{i\t})|^2\,
 \frac{d\t}{2\pi} \)^{1/2}<\infty\,.
$$
Any such function has radial limits, $f(e^{i\t}) =\lim_{r\to 1^-} f(r e^{i\t})$, almost everywhere on the unit circle $\T$ with respect to the normalized arc length measure: $dm(\t)=d\t/(2\pi)$, and its norm can be recovered from the boundary values, as well as from the power series $f(z)= \sum_{n=0}^\infty a_n z^n$ in the disk, as follows:
\begin{equation}
 \|f\|_2 = \( \int_{\T} |f|^2\,d m \)^{1/2} = \( \sum_{n=0}^\infty |a_n|^2 \)^{1/2}\,.
 \label{h2-norm}
\end{equation}
A function $\f$ defined in $\D$ is said to be a \textit{pointwise multiplier\/} of $H^2$ into itself if $\f f\in H^2$ for all $f\in H^2$; in other words, if the \textit{multiplication operator\/} $M_\f$, defined by $M_\f (f)=\f f$, maps $H^2$ into itself. By standard pointwise estimates for the functions in $H^2$ and the Closed Graph Theorem, this automatically implies boundedness of the operator $M_\f$. It is also known that this happens if and only if $\f\in H^\infty$ and, furthermore, $M_\f=\|\f\|_\infty$. For the discussions of this, we refer the reader to \cite{ADMV} or \cite{V}, for example.
\par
\subsection*{A characterization of inner functions in terms of their Taylor coefficients}
As observed above, a function $\f$ analytic in $\D$ is bounded by one if and only if the multiplication operator $M_\f$ has norm at most one on $H^2$. This was the key to Kortram's short proof \cite{K} of Schur's criterion. This idea can also be adapted to the context of inner functions - with equalities and requiring slightly modified sums.
\begin{thm} \label{thm-schur}
Let $\f$ be analytic in the unit disk, $\f(z)=\sum_{n=0}^\infty \la_n z^n$. Then the following conditions are equivalent:
\begin{itemize}
\item[(a)]
$\f$ is an inner function;
\item[(b)]
$M_\f$ is an isometric pointwise multiplier of $H^2$ into itself; in other words, $\|\f f\|_2 = \|f\|_2$ for all $f\in H^2$;
\item[(c)]
$$
 \sum_{n=0}^\infty \left| \sum_{j=0}^n a_j \lambda_{n-j} \right|^2 = \sum_{n=0}^\infty |a_n|^2
$$
holds for all square-summable sequences $(a_n)_{n=1}^\infty$.
\item[(d)]
$$
 \sum_{n=0}^N \left| \sum_{j=0}^n a_j \lambda_{n-j} \right|^2 + \sum_{n=N+1}^\infty \left| \sum_{j=0}^N a_j \lambda_{n-j} \right|^2 = \sum_{n=0}^N |a_n|^2
$$
holds for all non-negative integers $N$ and for all possible choices $(a_0,a_1,\ldots,a_N) \in \C^{N+1}$.
\end{itemize}
\end{thm}
\par
That (a) is equivalent to (b) is contained in Theorem~2.1 from \cite{ADMV}. Here we give an indirect proof of this fact. The density of the polynomials in $H^2$ is in the background of the proof.
\par
\begin{proof}
It suffices to show that (a) $\Rightarrow$ (b) $\Rightarrow$ (c) $\Rightarrow$ (d) $\Rightarrow$ (a).
\par
\fbox{(a) $\Rightarrow$ (b)} Since $|\f|=1$ almost everywhere on $\T$, for every $f\in H^2$ we have
$$
 \|\f f\|_2^2 = \int_\T |\f f|^2 dm = \int_\T |f|^2 dm = \|f\|_2^2\,.
$$
\par\smallskip
\fbox{(b) $\Rightarrow$ (c)} Let $f(z)= \sum_{n=0}^\infty a_n z^n$, where $\|f\|_2^2 = \sum_{n=0}^\infty |a_n|^2 < \infty$. Just note that the $n$-th coefficient of $\f f$ is none other than $\sum_{j=0}^n a_j \lambda_{n-j}$ by standard multiplication of two power series. The statement now follows from the assumption (b) and norm formula \eqref{h2-norm}.
\par\smallskip
\fbox{(c) $\Rightarrow$ (d)} This follows directly by substituting the stationary sequence $(a_0,a_1,\ldots,a_N,0,0,\ldots)$ into the equality in (c).
\par\smallskip
\fbox{(d) $\Rightarrow$ (a)} Proving this statement is the crux of the proof. To this end, first substitute $(a_0,a_1,\ldots,a_{N-1},a_N) =(0,0,\ldots,0,1)$ in (d) to obtain
$$
 \|\f\|_2^2 = |\lambda_0 |^2 + \sum_{n=N+1}^\infty | \lambda_{n-N} |^2 = |a_N|^2 = 1\,.
$$
Now let $f$ be an arbitrary function in $H^2$ with $f(z)= \sum_{n=0}^\infty a_n z^n$ in the disk. Apply equality (d) to the $(N+1)$-tuple  $(a_0,a_1,\ldots,a_N)$ of initial coefficients of $f$ to deduce that
$$
 \sum_{n=0}^N \left| \sum_{j=0}^n a_j \lambda_{n-j} \right|^2 \le  \sum_{n=0}^N |a_n|^2\,.
$$
Do this for every $N$ and then let $N\to\infty$. By \eqref{h2-norm}, we get $\|\f f\|_2\le \|f\|_2$. This means that the multiplication operator $M_\f$ has norm at most one and, thus, $\|\f\|_\infty\le 1$.
Hence in the obvious double inequality:
$$
 1 = \|\f\|_2 \le \|\f\|_\infty \le 1
$$
(which follows from \eqref{h2-norm} and our assumptions on $\f$)
equality must hold throughout, so $|\f|$ must have constant value one almost everywhere on the unit circle, which proves (a).
\end{proof}
\par
\subsection*{A sequential interpretation of Beurling's theorem}
In view of Beurling's theorem, for every invariant subspace $M$ of $H^2$ there is an inner function $\f$ such that for all $g\in M$ we have $g=\f f$, and clearly $\f\in M$. Writing the power series expansions and using the above description of inner functions, we can immediately translate the statement of Theorem~\ref{thm-schur} into the language of $\ell^2$ sequences.
\begin{cor} \label{cor-beur}
A closed non-trivial invariant subspace $M$ of $\ell^2$ is invariant for the shift operator if and only if there exists a fixed sequence $(\la_n)_{n=0}^\infty$ in $M$ which satisfies either of the conditions (c) or (d) of Theorem~\ref{thm-schur} (hence, both of them) and such that every sequence in $M$ is of the form
$$
 \( \sum_{k=0}^n \la_k c_{n-k} \)_{n=0}^\infty
$$
for some $(c_n)_{n=0}^\infty \in\ell^2$.
\end{cor}
\par



\begin{thebibliography}{ADMV}

\bibitem{AJ}
P. Ahern and M. Jevti\'c, Mean modulus and the fractional derivative of an inner function, \textit{Complex Variables Theory Appl.\/} \textbf{3}  (1984), no. 4, 431–445.

\bibitem{ADMV}
A. Aleman, P. Duren, M.J. Mart\'{\i}n, and D. Vukoti\'c, Multiplicative Isometries and Isometric Zero-Divisors, \textit{Canad.J. Math.\/} \textbf{62} (2010), No. 5, 961–974.

\bibitem{B}
A. Beurling, On two problems concerning linear transformations on Hilbert space, \textit{Acta Math.} \textbf{81} (1949), 239--255.

\bibitem{D}
P.L. Duren, \textit{Theory of $H\sp p$ Spaces}, Academic Press, New~York 1970; reprinted by Dover, Mineola, NY 2000.

\bibitem{G}
J.B. Garnett, \emph{Bounded Analytic Functions}, Academic Press, New York 1981.

\bibitem{K}
R.A. Kortram, A simple proof for Schur's theorem, \textit{Proc. Amer. Math. Soc.} \textbf{129} (2001), No. 11, 3211--3212.

\bibitem{NS}
D.J.  Newman and H.S. Shapiro, The Taylor coefficients of inner functions,  \textit{Michigan Math. J.\/} \textbf{9} (1962), 249–255.

\bibitem{S}
I. Schur, \"Uber Potenzreihen die im Innern des Einheitskreises beschr\"ankt sind, \textit{J. Reine Angew. Math.\/} \textbf{147} (1917), 205-232.

\bibitem{T}
M. Tsuji, \textit{Potential Theory in Modern Function Theory}, Maruzen, Tokyo 1959.

\bibitem{V}
D. Vukoti\'c, Analytic Toeplitz operators on the Hardy space $H^p$ - a survey,  \textit{Bull. Belgian Math. Soc. -- Simon Stevin\/}, \textbf{10}  (2003), 101--113.

\end{thebibliography}
\end{document}